\documentclass[11pt]{article}
\usepackage{amsmath}
\usepackage{amssymb}
\usepackage{graphics}
\usepackage{graphicx}
\allowdisplaybreaks[4]
\newtheorem{theo}{\hspace*{\parindent}Theorem}

\newtheorem{lemma}{\hspace*{\parindent}Lemma}

\def\cap{\mathrm{cap}}
\def\Arg{\mathrm{Arg}}
\def\C{\overline{\mathbb{C}}}
\def\Z{\mathbb{Z}}
\def\R{\mathbf{R}}
\def\Z{\mathrm{Z}}
\def\vt{\vartheta}
\def\ve{\varepsilon}
\def\O{{\cal O}}
\def\sn{\mathrm{sn}}
\def\cn{\mathop{\mathrm{cn}}}
\def\dn{\mathop{\mathrm{dn}}}
\def\union{\mathop{\cup}}
\def\cp{\mathrm{Cap}}
\def\adm{\mathrm{adm}}

\title{Capacity of a condenser whose plates are circular arcs}
\author{D. Karp\footnote{Institute of Applied Mathematics, Far Eastern Branch of
The Russian Academy of Sciences, 7 Radio Street, Vladivostok,
690041, Russia. e-mail:\emph{dmkrp@yandex.ru}}}
\date{}
\begin{document}
\maketitle

\begin{center}
\parbox{12cm}{
\small\textbf{Abstract.} We find an asymptotic formula for the
conformal capacity of a plane condenser both plate of which are
concentric circular arcs as the distance between them vanishes.
This result generalizes the formula for the capacity of parallel
linear plate condenser found by Simonenko and Chekulaeva in 1972
and sheds light on the problem of finding an asymptotic formula
for the capacity of condenser whose plates are arbitrary parallel
curves. This problem was posed and partially solved by R.
K\"{u}hnau in 1998.}
\end{center}


\bigskip
\paragraph{1.~Introduction.}
A pair of closed non-empty and non-intersecting subsets $E_0$,
$E_1$ of the extended complex plane $\C$ will be called {\em a
condenser} and denoted by $C=(E_0,E_1)$. The sets $E_0$, $E_1$ are
called {\em the plates} and the open set
$G=\C\backslash(E_0\union{E_1})$ is called {\em the field} of the
condenser $C$.  \emph{Capacity} of $C$ is defined by
\begin{equation}\label{eq:Cap-defined}
\cp(C)=\inf\limits_{u\in\adm(C)}\int\limits_{\mathbb{C}}|\nabla{u}|^2d\sigma(z),
\end{equation}
where integration is with respect to the planar Lebesgue measure,
and $\adm(C)$ denotes the collection of continuous functions
$\C\to\R$ satisfying Lipschitz condition in a neighborhood of
every finite point in $G$, possibly excluding a finite number of
points, and such that $u(z)=i$ for $z\in{E_i}$, $i=0,1$.  The
function $\omega$ which is harmonic in $G$ and assumes the value
$i$ on $E_i$ is called the {\em potential function} of the
condenser $C$. According to the Dirichlet principle it solves the
extremal problem (\ref{eq:Cap-defined}).  The potential function
exists for all condensers encountered in this paper.  The lines
orthogonal to the level curves of $\omega$ (or equivalently the
level curves of the harmonic conjugate of $\omega$) are called
\emph{the field lines} of the condenser $C$. The crucial property
of the family $\Gamma$ comprising the field lines of $C$ is that
its modulus $M(\Gamma)$ (the number reciprocal to the extremal
length, see \cite{Ahlfors}) is equal to the capacity of $C$
\cite{Gehring,Ziemer}:
\begin{equation}\label{eq:modulus-capacity}
M(\Gamma)=\cp(C).
\end{equation}
Both quantities in the above equality are conformal invariants.
For the purposes of this paper it is sufficient to mention that
the modulus of the family of straight line segments connecting the
sides of length $a$ of a rectangle and going parallel to the side
of length $b$ equals $a/b$ \cite{Ahlfors}.

Condensers defined above and their various generalizations  found
numerous applications in geometric theory of functions and many
other areas. The study of the asymptotic behavior of the capacity
of a specifically designed condenser led V.N.~Dubinin to the
definition and calculation of the generalized reduced modulus
particular cases of which go back to Gr\"{o}tzsch,
Teichm\"{u}ller, Ahlfors and Beurling. Using this concept Dubinin
and his students proved a number of covering and distortion
theorems for analytic functions and gave solutions to several
previously unsolved extremal partition problems. Many new
inequalities for polynomials and rational functions have been also
established. See for instance \cite{Dubinin} and references
therein.

The condenser $C_0=(E_0,E_1)$ whose plates are parallel linear
segments
\[
E_0=\{z:z\in[-L/2,L/2]\},~~E_1=\{z:z\in[-L/2+ih,L/2+ih]\},
\]
was studied by Simonenko and Chekulaeva in \cite{Simonenko}. This
condenser may also be thought of as a three dimensional condenser
consisting of two parallel infinite bands. The authors used the
following transcendental equation for the capacity of $C_0$ (which
is essentially contained already in \cite{Betz}):
\begin{equation}\label{eq:SC-equation}
\frac{\pi L}{2h}=KE(\phi,k)-EF(\phi,k),
\end{equation}
where
\[
\phi=\arcsin\left(\frac{1}{k}\sqrt{1-E/K}\right)~~\text{and}~~\cp(C_0)=K/K'.
\]
Here $K=K(k)$ is the complete elliptic integrals of the first kind
(see (\ref{eq:K-defined})), $K'=K(\sqrt{1-k^2})$, $E=E(k)$ is the
complete elliptic integral of the second kind and $F(\phi,k)$,
$E(\phi,k)$ are the incomplete elliptic integrals of the first and
second kind, respectively (see (\ref{eq:F-defined}),
(\ref{eq:E-defined})). From this equation Simonenko and Chekulaeva
derived an asymptotic formula for the capacity of $C_0$ when
plates approach each other (i.e. $h\to{0}$) which in our notation
can be written as:
\[
\cp(C_0)=\frac{L}{h}+\frac{1}{\pi}\ln\frac{1}{h}+\frac{1}{\pi}\left(1+\ln(2\pi{L})\right)
+\frac{h}{\pi^2L}\ln\frac{1}{h}+\frac{h}{\pi^2L}\left(\frac{1}{2}+\ln(2\pi{L})\right)
\]
\begin{equation}\label{eq:SC-asymp}
 -\frac{h^2}{2\pi^3L^2}\ln^2\frac{1}{h}
+\frac{h^2}{\pi^3L^2}\left(\frac{1}{2}-\ln(2\pi{L})\right)\ln\frac{1}{h}+\O(h^2).
\end{equation}
Here and below $f(x)=\O(g(x))$ as $x\to{a}$ means that
$C_1g(x)\leq|f(x)|\leq{C_2}g(x)$ for all $x$ in a neighbourhood of
$a$.

Reiner K\"{u}hnau considered in \cite{Kuehnau2} a general parallel
plate condenser $C_K=(E_0,E_1)$ whose plates are arbitrary
parallel curves.  More precisely, let a Jordan curve parameterized
by the arc length $s\in[0,L]$ be defined by a three times
differentiable function $z=z(s)$. The plates $E_0$, $E_1$ of the
condenser $C_K$ are curves parallel to $z(s)$ and located at equal
distances from it:
\[
E_0=\{z:z_+(s)=z(s)+ihz'(s)/2,~s\in[0,L]\},
~~~E_1=\{z:z_-(s)=z(s)-ihz'(s)/2,~s\in[0,L]\}.
\]
Using a properly chosen family of curves K\"{u}hnau has proved the
asymptotic formula
\begin{equation}\label{eq:C-K-asymp}
\cp(C_K)=\frac{L}{h}+\frac{1}{\pi}\ln\frac{1}{h}+\O(1).
\end{equation}

The goal of this paper is to find an asymptotic expansion for the
capacity of the condenser $C_\rho$ which is a particular case of
K\"{u}hnau's condenser when $z(s)$ traverses a circular arc
symmetric with respect to real axis:
\[
z(s)=\rho\exp[i(s/\rho-\gamma)],~~s\in[0,L],~~L=2\gamma\rho.
\]
Differentiation yields:
\begin{equation}\label{eq:zpm}
z_{\pm}(s)=(\rho{\mp}h/2)\exp[i(s/\rho-\gamma)].
\end{equation}
Hence the plates of the condenser $C_\rho=(E_0,E_1)$ are
\begin{equation}\label{eq:Crho-defined}
E_0=\{z:z=z_+(s),s\in[0,L]\},~~~E_1=\{z:z=z_-(s),s\in[0,L]\}.
\end{equation}
K\"{u}hnau's approximation for its capacity is given by
(\ref{eq:C-K-asymp}).  The main result of this paper is
\begin{theo}\label{th:main}
For $h\to{0}$ the following asymptotic expansion holds true
\[
\cp(C_\rho)=\frac{L}{h}+\frac{1}{\pi}\ln\frac{1}{h}+\frac{1}{\pi}\left(1+\ln(4\pi\rho\sin(L/2\rho))\right)+
\frac{\cot(L/2\rho)}{2\pi^2\rho}h\ln(1/h)
+\frac{\cot(L/2\rho)}{2\pi^2\rho}h\times
\]
\begin{equation}\label{eq:C-rho-asymp}
\times\!\!\left[\frac{1}{2}+\ln(4\pi\rho\sin(L/2\rho))\right]-\frac{h^2\ln^2(1/h)}{8\pi^3\rho^2\sin^2(L/2\rho)}
-\frac{h^2\ln(1/h)\left(2\ln(4\pi\rho\sin(L/2\rho))-\cos(L/\rho)\right)}{8\pi^3\rho^2\sin^2(L/2\rho)}+\O(h^2).
\end{equation}
\end{theo}
Our method allows one to compute as many terms in the above
asymptotic expansion as one wishes. Formula (\ref{eq:SC-asymp})
for the capacity of parallel linear plate condenser follows from
(\ref{eq:C-rho-asymp}) if we let $\rho\to\infty$. We have also
obtained an analogue of equation (\ref{eq:SC-equation}) - see
Theorem~\ref{th:cap-equation}. Formula (\ref{eq:gamma*-found})
from this Theorem has been announced in \cite{Dymch}.  An equation
for the capacity of a condenser one plate of which is a circular
arc and the other is a radius of the same circle was found in
\cite{Davy}.

Related developments in three dimensions were considered in
\cite{Kuehnau1} and  \cite{Soibel}.  The latter paper gives an
asymptotic expansion for the capacity of a condenser whose plates
are parallel planar figures of arbitrary shape as the distance
between plates vanishes.  Specific calculation of coefficients,
however, has only been done in the classical case of the circular
plates. In \cite{Kuehnau1} R.~K\"{u}hnau gives an asymptotic
formula for the capacity of a condenser of which both plates are
parallel finite surfaces when its field is restricted to the space
between them.

\begin{figure}
\centering \scalebox{1}{\includegraphics{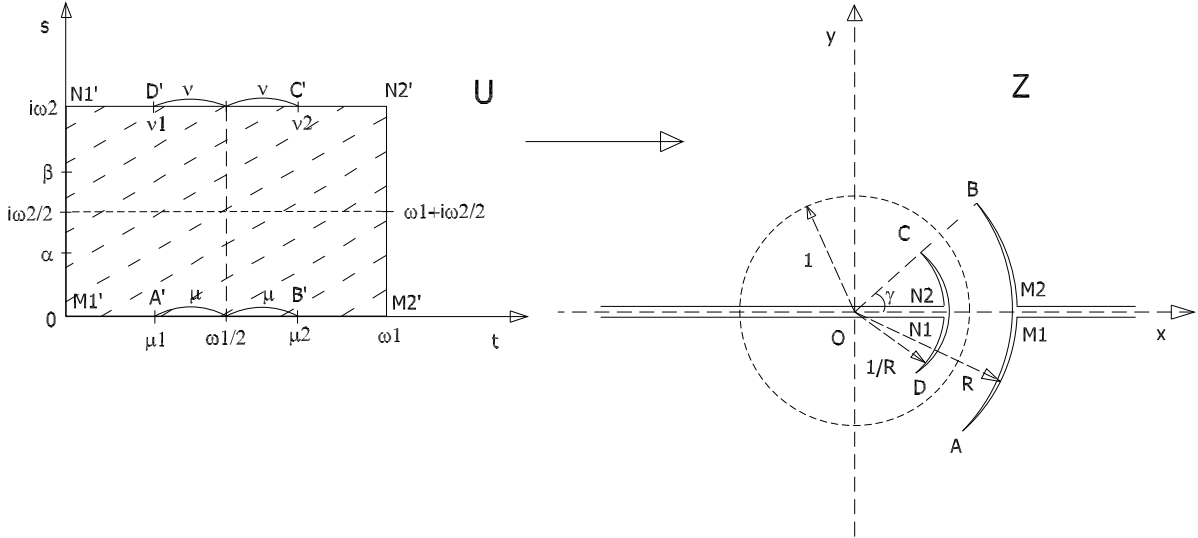}}
\caption{The mapping realized by $z(u)$} \label{fig:z-u}
\end{figure}

\paragraph{2. Preliminaries.}  To derive a formula for the
capacity of $C_\rho$ it is sufficient to find one for $C_1$.  Once
this has been done the conformal mapping $z\to\rho{z}$ will bring
the result for $C_\rho$.  Thus we take $\rho=1$, $L=2\gamma$ and
$\delta=h/2$. Before proceeding further we transform our condenser
$C_1$ into the condenser $C_1'$ which will be more convenient to
work with. To this end carry out the conformal mapping
\[
z\to\frac{z}{\sqrt{1-\delta^2}}.
\]
Under this mapping the condenser $C_1$ transforms into condenser
$C_1'$ with plates
\begin{equation}\label{eq:E0-C1prime}
E_0=\left\{z:z=e^{i\phi}/R,~R=\sqrt{\frac{1+\delta}{1-\delta}}~,\phi\in[-\gamma,\gamma]\right\},
\end{equation}
\begin{equation}\label{eq:E1-C1prime}
E_1=\left\{z:z=Re^{i\phi},~R=\sqrt{\frac{1+\delta}{1-\delta}}~,\phi\in[-\gamma,\gamma]\right\}
\end{equation}
(the new variable has been again denoted by $z$). Put $\ve=R-1$.
Then:
\begin{equation}\label{eq:eps-delta}
\ve=R-1=\sqrt{\frac{1+\delta}{1-\delta}}-1=\delta\left(1+\frac{\delta}{2}+\frac{\delta^2}{2}\right)+\O(\delta^4)
\end{equation}
and
\begin{equation}\label{eq:delta-eps}
\delta=\frac{(1+\ve)^2-1}{(1+\ve)^2+1}=\ve\left(1-\frac{\ve}{2}+\frac{\ve^3}{4}\right)+\O(\ve^5).
\end{equation}
Conformal invariance of capacity combined with
(\ref{eq:C-K-asymp}) and (\ref{eq:delta-eps}) yield
($\delta=h/2$):
\begin{equation}\label{eq:C1-first-asymp}
\cp(C_1)=\cp(C_1')=\frac{\gamma}{\ve}+\frac{1}{\pi}\ln\frac{1}{\ve}+\O(1).
\end{equation}
To find out more about the capacity of $C_1'$ we map the interior
of the rectangle with vertices $0$, $\omega_1$,
$\omega_1+i\omega_2$, $i\omega_2$ in the complex $\C_u$ plane onto
the field of $C_1'$ endowed with a slit along the real axis
connecting the plates and going through infinity (see
Figure~\ref{fig:z-u}). The modulus $M(\Gamma)$ of the family
$\Gamma$ comprising vertical line segments which connect the top
and the bottom side of the rectangle equals $\omega_1/\omega_2$.
The conformal image $\Gamma'$ of $\Gamma$ under the mapping
$u\to{z(u)}$ is the family of field lines of the condenser $C_1'$.
Obviously, the slit cannot change the modulus of $\Gamma'$ since
no field line can cross it due to symmetry. Conformal invariance
and identity (\ref{eq:modulus-capacity}) give:
\begin{equation}\label{eq:capacity-basic}
\cp(C_1')=M(\Gamma')=M(\Gamma)=\omega_1/\omega_2.
\end{equation}
The mapping $z(u)$ will be constructed in the lemma below. Before
we proceed to formulating it let us remind the definitions of
theta-functions \cite{Akhiezer,Bat3}:
\begin{equation}\label{eq:theta1-def}
\vt_1(w;q)=2\sum\limits_{n=0}^{\infty}(-1)^nq^{(n+1/2)^2}\sin((2n+1)w)=
2q^{1/4}\sin{w}\prod\limits_{n=1}^{\infty}(1-q^{2n})(1-2q^{2n}\cos{2w}+q^{4n}),
\end{equation}
\begin{equation}\label{eq:theta2-def}
\vt_2(w;q)=2\sum\limits_{n=0}^{\infty}q^{(n+1/2)^2}\cos((2n+1)w)=
2q^{1/4}\cos{w}\prod\limits_{n=1}^{\infty}(1-q^{2n})(1+2q^{2n}\cos{2w}+q^{4n}),
\end{equation}
\begin{equation}\label{eq:theta3-def}
\vt_3(w;q)=1+2\sum\limits_{n=1}^{\infty}q^{n^2}\cos(2nw)=
\prod\limits_{n=1}^{\infty}(1-q^{2n})(1+q^{2n-1}e^{2iw})(1+q^{2n-1}e^{-2iw}),
\end{equation}
\begin{equation}\label{eq:theta4-def}
\vt_4(w;q)=1+2\sum\limits_{n=1}^{\infty}(-1)^nq^{n^2}\cos(2nw)=
\prod\limits_{n=1}^{\infty}(1-q^{2n})(1-q^{2n-1}e^{2iw})(1-q^{2n-1}e^{-2iw}).
\end{equation}
The following properties of $\vt_i$ will prove useful in the
sequel \cite{Akhiezer,Bat3}(we will omit $q$ in $\vt_i(w;q)$ when
it cannot lead to confusion):
\begin{equation}\label{eq:theta-prop1}
\vt_4(w+\pi)=\vt_4(w)
\end{equation}
\begin{equation}\label{eq:theta-prop2}
\vt_4(w+\pi\tau)=-q^{-1}e^{-2iw}\vt_4(w),~~~~\tau=\frac{1}{i\pi}\ln{q},
\end{equation}
\begin{equation}\label{eq:theta-prop3}
\vt_4(w+\pi\tau/2;q)=ie^{-i\pi\tau/4}e^{-iy}\vt_1(w;q),
\end{equation}
\begin{equation}\label{eq:theta-prop4}
\vt_4(w+\pi/2;q)=\vt_3(w;q),
\end{equation}
\begin{equation}\label{eq:theta-prop5}
\vt_i(\overline{w})=\overline{\vt_i(w)},~~~~i=1,3,4,
\end{equation}
\begin{equation}\label{eq:theta-prop6}
\Im(\vt_i(ix))=0,~~~~x\in\R,~~~~i=1,3,4,
\end{equation}
\begin{equation}\label{eq:theta-prop7}
\vt_i(-w))=\vt_i(w),~~~i=3,4,
\end{equation}
\begin{equation}\label{eq:theta-prop8}
\vt_1(-w))=-\vt_1(w),
\end{equation}
\begin{equation}\label{eq:theta-prop9}
\vt_4(\pi\tau/2+\pi{m}+\pi{n}))=0,~~~~~~~n,m\in\Z.
\end{equation}
\begin{lemma}\label{lemma:mapping}
Let $\omega_1>0$, $\omega_2>0$ and $1<R<e^{\pi\omega_2/\omega_1}$.
The function
\begin{equation}\label{eq:z-defined}
z(u)=R\frac{\vt_4({\pi}u/\omega_1-(i\ln{R})/2;q)}{\vt_4({\pi}u/\omega_1+(i\ln{R})/2;q)},
~~~q=e^{i\pi\tau},~~~\tau=i\frac{\omega_2}{\omega_1},
\end{equation}
maps the rectangle with vertices $0$, $\omega_1$,
$\omega_1+i\omega_2$, $i\omega_2$ in the complex $\C_u$ plane
conformally and univalently onto the entire $\C_z$ plane with two
circular slits with radii $R$ and $1/R$ symmetric with respect to
real axis and connected by the slit along the real axis going
through infinity.
\end{lemma}
The mapping (\ref{eq:z-defined}) was essentially constructed in
\cite{Sedov}.
For completeness we give a direct proof  here. \\
\textbf{Proof.} Let $u$ traverse the boundary of the rectangle
shown on Figure~\ref{fig:z-u}.  To trace its image $z(u)$ we first
note that by (\ref{eq:theta-prop1}) and (\ref{eq:theta-prop2})
($\tau=i\omega_2/\omega_1$):
\begin{equation}\label{eq:z-periodic}
z(u+\omega_1)=z(u),
\end{equation}
\begin{equation}\label{eq:z-omega2}
z(u+\omega_2)=R\frac{\vt_4({\pi}u/\omega_1-(i\ln{R})/2+\pi\tau;q)}{\vt_4({\pi}u/\omega_1+(i\ln{R})/2+\pi\tau;q)}=
\frac{-q^{-1}e^{-2i(\pi{u}/\omega_1-(i\ln{R})/2)}}{-q^{-1}e^{-2i(\pi{u}/\omega_1+(i\ln{R})/2)}}z(u)=\frac{z(u)}{R^2}.
\end{equation}
For the bottom side $(M_1',M_2')$ write $u=t\in[0,\omega_1]$. We
employ (\ref{eq:theta-prop3}) to get:
\[
|z(t)|=R\left|\frac{\vt_4({\pi}t/\omega_1-(i\ln{R})/2)}{\vt_4({\pi}t/\omega_1+(i\ln{R})/2)}\right|=R
\]
and
\[
z(0)=z(\omega_1)=R\frac{\vt_4(-i\ln{R}/2)}{\vt_4(i\ln{R}/2)}=R\frac{\overline{\vt_4(i\ln{R}/2)}}{\vt_4(i\ln{R}/2)}=R
\]
by (\ref{eq:theta-prop5}) and (\ref{eq:theta-prop6}).  It follows
that the interval $[0,\omega_1]$ is mapped onto a circular arc
beginning and ending at the point $z=R$. This arc is symmetric
with respect to real axis since for any
$\lambda\in(0,\omega_1/2)$:
\[
z(\omega_1/2+\lambda)=R\frac{\vt_4(\pi\lambda/\omega_1-(i\ln{R})/2)+\pi/2)}{\vt_4(\pi\lambda/\omega_1+(i\ln{R})/2)+\pi/2)}=
R\frac{\vt_3(\pi\lambda/\omega_1-(i\ln{R})/2))}{\vt_3(\pi\lambda/\omega_1+(i\ln{R})/2)}=
\]\[
=R\frac{\overline{\vt_3(\pi\lambda/\omega_1+(i\ln{R})/2))}}{\overline{\vt_3(\pi\lambda/\omega_1-(i\ln{R})/2)}}=
R\frac{\overline{\vt_3(-\pi\lambda/\omega_1-(i\ln{R})/2))}}{\overline{\vt_3(-\pi\lambda/\omega_1+(i\ln{R})/2)}}
=\overline{z(\omega_1/2-\lambda)}
\]
by (\ref{eq:theta-prop5}) and (\ref{eq:theta-prop7}). In
particular, $z(\omega_1/2)=R$.

For the top side $(N_1',N_2')$ write $u=t+i\omega_2$,
$t\in[0,\omega_1]$, and  by (\ref{eq:z-omega2}):
\[
|z(t+i\omega_2)|=\frac{|z(t)|}{R^2}=\frac{1}{R}.
\]
Previous calculations for the bottom side combined with
(\ref{eq:z-omega2}) show that the side $(N_1',N_2')$ is mapped
onto a circular arc beginning and ending at the point $z=1/R$
which is symmetric with respect to real axis.

When $u\in[0,i\omega_2]$ we see that $z(u)$ is real by
(\ref{eq:theta-prop6}).  Due to periodicity (\ref{eq:z-periodic})
the same values are assumed by $z(u)$ when
$u\in[\omega_1,\omega_1+i\omega_2]$.  Suppose $\alpha$ is the
preimage of infinity, so that $z(\alpha)=\infty$, and $\beta$ is
preimage of origin, so that $z(\beta)=0$. Then according to
(\ref{eq:theta-prop9}):
\[
\frac{\pi\alpha}{\omega_1}+\frac{1}{2}i\ln{R}=\frac{i\pi\omega_2}{2\omega_1}~~\Leftrightarrow~~
\alpha=i\left(\frac{\omega_2}{2}-\frac{\omega_1}{2\pi}\ln{R}\right)
\]
and
\[
\frac{\pi\beta}{\omega_1}-\frac{1}{2}i\ln{R}=\frac{i\pi\omega_2}{2\omega_1}~~\Leftrightarrow~~
\beta=i\left(\frac{\omega_2}{2}+\frac{\omega_1}{2\pi}\ln{R}\right).
\]
Hence $\alpha$ lies in $(0,i\omega_2/2)$ when
$1<R<e^{\pi\omega_2/\omega_1}$ and $\alpha$ and $\beta$ are
symmetric with respect to $i\omega_2/2$.  For $u=t+i\omega_2/2$,
$t\in[0,\omega_1]$, we get by (\ref{eq:theta-prop3}):
\begin{multline*}
z(t+i\omega_2/2)=R\frac{\vt_4(\pi{t}/\omega_1-(i\ln{R})/2+\pi\tau/2)}{\vt_4(\pi{t}/\omega_1+(i\ln{R})/2+\pi\tau/2)}\\
=R\frac{ie^{i\pi\tau/4-i(\pi\tau/\omega_1-(i\ln{R})/2)}}{ie^{i\pi\tau/4-i(\pi\tau/\omega_1+(i\ln{R})/2)}}
\frac{\vt_1(\pi{t}/\omega_1-(i\ln{R})/2)}{\vt_1(\pi{t}/\omega_1+(i\ln{R})/2)}=
\frac{\vt_1(\pi{t}/\omega_1-(i\ln{R})/2)}{\vt_1(\pi{t}/\omega_1+(i\ln{R})/2)}.
\end{multline*}
Hence by (\ref{eq:theta-prop5})
\[
|z(t+i\omega_2/2)|=1
\]
and by (\ref{eq:theta-prop8})
\[
z(i\omega_2/2)=z(\omega_1+i\omega_2/2)=-1.
\]
Thus the dotted line on Figure~\ref{fig:z-u} connecting
$i\omega_2/2$ and $\omega_1+i\omega_2/2$ is mapped onto the unit
circle in $\C_z$ plane.

Finally, we see that both arcs have the same angular spread, so
that the points $C$ and $B$ lie on the same beam.  This is an
obvious consequence of (\ref{eq:z-omega2}) since
\[
\arg[z(t+i\omega_2)]=\arg[z(t)/R^2]=\arg[z(t)].
\]
It follows that  $\nu=\mu$ on Figure~\ref{fig:z-u}. We denote the
angle between real axis and the beam $OCB$ by $\gamma$ so that
\begin{equation}\label{eq:gamma*}
\gamma\equiv\max\limits_{t\in(0,\omega1)}\arg[z(t)].
\end{equation}
Above speculations are summarized the in the following chart:
\[
\renewcommand\arraystretch{1.2}
\begin{array}{|c|c|c|c|c|c|c|c|c|c|}
\hline &&&&&&&&&\\[-13pt]
u &  0 & \mu_1 & \omega_1/2 & \mu_2 & \omega_1
&\omega_1+\alpha&\omega_1+i\omega_2/2 & \omega_1+\beta &
\omega_1+i\omega_2 \\[2pt]
\hline z & R & Re^{-i\gamma}  & R &  Re^{i\gamma}& R &\infty
&-1 & 0 &1/R\\[2pt]
\hline
\end{array}
\]
\[
\renewcommand\arraystretch{1.2}
\begin{array}{|c|c|c|c|c|c|c|c|}
\hline &&&&&&&\\[-13pt]
u &\nu_2 &\omega_1/2+i\omega_2& \nu_1 &i\omega_2 &\beta& i\omega_2/2 &\alpha \\[2pt]
\hline z & e^{i\gamma}/R & 1/R &e^{-i\gamma}/R &1/R &0 &-1& \infty\\[2pt]
\hline
\end{array}
\]
This completes the proof of the lemma.

\paragraph{3.~Equation for the capacity of $C_1'$.}  In this section we will
derive a transcendental equation for $\cp(C_1')$ which can be
viewed as a generalization of (\ref{eq:SC-equation}). We will
regard $\gamma$ and $R$ (and hence $\cp(C_1')$) as being fixed.
From (\ref{eq:capacity-basic}) and (\ref{eq:z-defined}) $q$ is
also fixed and computed by
\begin{equation}\label{eq:q}
q=e^{-\pi/\cap(C_1')}.
\end{equation}\
We still have a free scaling factor $\omega_1$. Choose
\begin{equation}\label{eq:omega1-fixed}
\omega_1=2K(q)=\pi\vt_3^2(0;q)=\pi\left(1+2\sum\limits_{n=1}^{\infty}q^{n^2}\right)^2.
\end{equation}
Then according to (\ref{eq:capacity-basic}) and (\ref{eq:q}):
\begin{equation}\label{eq:omega2-fixed}
\omega_2=\frac{2}{\pi}K(q)\ln{\frac{1}{q}}=\vt_3^2(0;q)\ln{\frac{1}{q}}=\ln{\frac{1}{q}}\left(1+2\sum\limits_{n=1}^{\infty}q^{n^2}\right)^2=2K'(q).
\end{equation}
$K$ and $K'$ defined by (\ref{eq:omega1-fixed}) and
(\ref{eq:omega2-fixed}) are known to be the complete elliptic
integrals of the moduli
\begin{equation}\label{eq:k-defined}
k^2=\frac{\vt_2^4(0;q)}{\vt_3^4(0;q)}=16\left[\frac{\sum\limits_{n=0}^{\infty}q^{(n+1/2)^2}}{1+2\sum\limits_{n=1}^{\infty}q^{n^2}}\right]^4
\end{equation}
and
\begin{equation}\label{eq:k'-defined}
{k'}^2=1-k^2=\frac{\vt_4^4(0;q)}{\vt_3^4(0;q)}=\left[\frac{1+2\sum\limits_{n=1}^{\infty}(-1)^nq^{n^2}}{1+2\sum\limits_{n=1}^{\infty}q^{n^2}}\right]^4,
\end{equation}
respectively \cite[\S{30}]{Akhiezer}, that is:
\begin{equation}\label{eq:K-defined}
K=K(k)=\int\limits_{0}^{1}\frac{dt}{\sqrt{(1-t^2)(1-k^2t^2)}},~~~K'=K(k')=\int\limits_{0}^{1}\frac{dt}{\sqrt{(1-t^2)(1-{k'}^2t^2)}}.
\end{equation}
Define
\begin{equation}\label{eq:theta4-new}
\theta_4(w;q)=\vt_4(\pi{w}/2K;q).
\end{equation}
Then by (\ref{eq:z-defined}):
\begin{equation}\label{eq:z-new}
z(u)=R\frac{\theta_4(u-i\alpha;q)}{\theta_4(u+i\alpha;q)},
\end{equation}
where we introduced
\begin{equation}\label{eq:a-defiend}
\alpha=\frac{1}{\pi}K\ln{R}>0.
\end{equation}
According to \cite[\S{29}]{Akhiezer}:
\begin{equation}\label{eq:Pi-z}
\Pi(u,i\alpha)=\frac{1}{2}\ln\frac{\theta_4(u-i\alpha;q)}{\theta_4(u+i\alpha;q)}+u\Z(i\alpha)=\frac{1}{2}\ln(z(u)/R)+u\Z(i\alpha),
\end{equation}
where
\[
\Pi(u,i\alpha)=k^2\sn(i\alpha)\cn(i\alpha)\dn(i\alpha)\int\limits_{0}^{u}\frac{\sn^2(t)dt}{1-k^2\sn^2(i\alpha)\sn^2(t)}
\]
\begin{equation}\label{eq:specialPi}
=k^2b\sqrt{(1-b^2)(1-k^2b^2)}\int\limits_{0}^{x}\frac{t^2dt}{(1-k^2b^2t^2)\sqrt{(1-t^2)(1-k^2t^2)}}
\end{equation}
and
\begin{equation}\label{eq:b-defined}
x=\sn(u,k),~~~~b=\sn(i\alpha,k).
\end{equation}
The Jacobi $Z$-function is defined by
\begin{equation}\label{eq:Z-defined}
Z(i\alpha)=\frac{\theta_4'(i\alpha;q)}{\theta_4(i\alpha;q)}=E(\sn(i\alpha),k)-\frac{E(k)}{K(k)}i\alpha
\end{equation}
where
\begin{equation}\label{eq:E-defined}
E(b,k)=\int\limits_{0}^{b}\frac{\sqrt{1-k^2t^2dt}}{\sqrt{1-t^2}}=\int\limits_{0}^{i\alpha}[\dn(t,k)]^2dt,
\end{equation}
is the incomplete elliptic integral of the second kind. The
function $\Pi(u,i\alpha)$ can be expressed in terms of incomplete
elliptic integrals of the first kind
\begin{equation}\label{eq:F-defined}
F(x,k)=\int\limits_{0}^{x}\frac{dt}{\sqrt{(1-t^2)(1-k^2t^2)}}
\end{equation}
and of the third kind
\begin{equation}\label{eq:Pi-defined}
\Pi(x,\nu,k)=\int\limits_{0}^{x}\frac{dt}{(1+\nu{t^2})\sqrt{(1-t^2)(1-k^2t^2)}}
\end{equation}
by means of
\begin{equation}\label{eq:Pi-Pi}
\Pi(u,i\alpha)=-\frac{1}{b}\sqrt{(1-b^2)(1-k^2b^2)}\left[F(x,k)-\Pi(x,-k^2b^2,k)\right].
\end{equation}
For $u\in[0,K]$ put $\gamma(u)=\Arg(z(u))\in(-\pi,\pi]$ - the
principal value of argument $z$.  Since $|z(u)|=R$
\[
\ln{z(u)}=\ln{R}+i\gamma(u).
\]
On the other hand from (\ref{eq:Pi-z}):
\[
\ln{z(u)}=2\Pi(u,i\alpha)-2u\Z(i\alpha)+\ln{R}.
\]
Since $b$ defined in (\ref{eq:b-defined}) equals
$i\sn(\alpha,k')/\cn(\alpha,k')$ according to \cite{Akhiezer,Bat3}
and so is an imaginary number, we see from (\ref{eq:Z-defined}),
(\ref{eq:E-defined}) and (\ref{eq:Pi-Pi}) that for positive real
$u$ and $\alpha$ the values taken by $\Pi(u,i\alpha)$ and
$Z(u,i\alpha)$ lie on the imaginary axis. Hence
\[
\gamma(u)=\Im\ln{z(u)}=2(\Pi(u,i\alpha)-u\Z(i\alpha))/i.
\]
The single extremum of $\gamma(u)$ for $u\in[0,K]$ is a minimum at
$u=\mu_1$ (equal to $-\gamma$, see Figure~\ref{fig:z-u}). Hence we
can find $\mu_1$ by setting $d\gamma(u)/du$ to zero:
\[
i\frac{d\gamma(u)}{du}=\frac{2k^2\sn(i\alpha)\cn(i\alpha)\dn(i\alpha)\sn^2(u)}{1-k^2\sn^2(i\alpha)\sn^2(u)}-2\Z(i\alpha)=0,
\]
where we used (\ref{eq:specialPi}).  Finally,
\begin{equation}\label{eq:lambda-defined}
\lambda^2(R,q)\equiv\sn^{2}(\mu_1)=\frac{Z(i\alpha)}{k^2\sn(i\alpha)[\cn(i\alpha)\dn(i\alpha)+\sn(i\alpha)\Z(i\alpha)]}
\end{equation}
and
\begin{equation}\label{eq:mu1}
\mu_1=F\left(\frac{\sqrt{Z(i\alpha)/\sn(i\alpha)}}{k\sqrt{\cn(i\alpha)\dn(i\alpha)+\sn(i\alpha)\Z(i\alpha)}},k\right)=F(\lambda,k).
\end{equation}
The value of $\gamma$ defined by (\ref{eq:gamma*})  is then
revealed from:
\[
\gamma=2i(\Pi(\mu_1,i\alpha)-\mu_1\Z(i\alpha)),
\]
or by (\ref{eq:Pi-Pi}) and (\ref{eq:b-defined}):
\begin{equation}\label{eq:gamma*-found}
\gamma=2i\frac{\cn(i\alpha)\dn(i\alpha)}{\sn(i\alpha)}\left[\Pi(\lambda(R,q),-k^2\sn^2(i\alpha),k)-F(\lambda(R,q),k)\right]-2iF(\lambda(R,q),k)\Z(i\alpha).
\end{equation}
Thus we arrive at the following statement.
\begin{theo}\label{th:cap-equation}
Let the condenser $C_1'=(E_0,E_1)$ be  defined by
\emph{(\ref{eq:E0-C1prime})}, \emph{(\ref{eq:E1-C1prime})}. Then
its capacity $\cp(C_1')$ satisfies transcendental equation
\emph{(\ref{eq:gamma*-found})} with $\lambda$, $\alpha$,  $k$, $K$
and $q$ defined by \emph{(\ref{eq:lambda-defined})},
\emph{(\ref{eq:a-defiend})},
 \emph{(\ref{eq:k-defined})}, \emph{(\ref{eq:omega1-fixed})}
and \emph{(\ref{eq:q})}, respectively.
\end{theo}

\paragraph{4.~Proof of Theorem~\ref{th:main}.}
The crucial step of the proof is to obtain  an asymptotic
expansion for $\cp(C_1')$ as $R\to{1}$, while $\gamma$ remains
fixed. After this will have been accomplished carrying out the
conformal mappings inverse to those that led from $C_\rho$ to
$C_1'$ will complete the proof.

First we introduce the new nome $q_1$ which tends to zero when
$q\to{1}$:
\begin{equation}\label{eq:q1}
q_1=e^{-i\pi/\tau}=e^{\pi^2/\ln(q)}=e^{-\pi\cp(C_1')}.
\end{equation}
In terms of $q_1$ we obtain using the Jacobi transformations (and
bearing in mind that $\sqrt{-i\tau}=\sqrt{\pi}[\ln(1/q_1)]^{-1/2}$
by  (\ref{eq:q1})):
\begin{equation}\label{eq:theta1-newnome}
\vt_1(z,q)=\frac{i}{\sqrt{\pi}}[\ln(1/q_1)]^{1/2}e^{z^2\ln(q_1)/\pi^2}\vt_1(iz\ln(q_1)/\pi,q_1),
\end{equation}
\begin{equation}\label{eq:theta2-newnome}
\vt_2(z,q)=\frac{1}{\sqrt{\pi}}[\ln(1/q_1)]^{1/2}e^{z^2\ln(q_1)/\pi^2}\vt_4(iz\ln(q_1)/\pi,q_1),
\end{equation}
\begin{equation}\label{eq:theta3-newnome}
\vt_3(z,q)=\frac{1}{\sqrt{\pi}}[\ln(1/q_1)]^{1/2}e^{z^2\ln(q_1)/\pi^2}\vt_3(iz\ln(q_1)/\pi,q_1),
\end{equation}
\begin{equation}\label{eq:theta4-newnome}
\vt_4(z,q)=\frac{1}{\sqrt{\pi}}[\ln(1/q_1)]^{1/2}e^{z^2\ln(q_1)/\pi^2}\vt_2(iz\ln(q_1)/\pi,q_1).
\end{equation}

Denote
\begin{equation}\label{eq:eta-defined}
\eta=\frac{1}{2\pi}\ln(R)\ln(1/q_1)=\frac{1}{2}\ln(R)\cp(C_1').
\end{equation}
 For the the Jacobi elliptic functions we obtain
($\alpha=K\ln(R)/\pi$):
\begin{equation}\label{eq:sn-q1}
\sn(i\alpha,k)=\frac{1}{\sqrt{k}}\frac{\vt_1(i\alpha\pi/(2K),q)}{\vt_4(i\alpha\pi/(2K),q)}=
\frac{i}{\sqrt{k}}\frac{\vt_1(\eta,q_1)}{\vt_2(\eta,q_1)},
\end{equation}
\begin{equation}\label{eq:cn-q1}
\cn(i\alpha,k)=\frac{\sqrt{k'}}{\sqrt{k}}\frac{\vt_2(i\alpha\pi/(2K),q)}{\vt_4(i\alpha\pi/(2K),q)}=
\frac{\sqrt{k'}}{\sqrt{k}}\frac{\vt_4(\eta,q_1)}{\vt_2(\eta,q_1)},
\end{equation}
\begin{equation}\label{eq:dn-q1}
\dn(i\alpha,k)=\sqrt{k'}\frac{\vt_3(i\alpha\pi/(2K),q)}{\vt_4(i\alpha\pi/(2K),q)}=
\sqrt{k'}\frac{\vt_3(\eta,q_1)}{\vt_2(\eta,q_1)}
\end{equation}
and using (\ref{eq:theta4-new}) and (\ref{eq:Z-defined}):
\begin{equation}\label{eq:Z-q1}
\Z(i\alpha,k)=\frac{i\ln(q_1)}{2K}\left(\frac{1}{\pi}\ln(R)+\frac{\vt_2'(\eta,q_1)}{\vt_2(\eta,q_1)}\right).
\end{equation}

Setting $\ve=R-1$, we immediately derive from
\[
\ln(1+\ve)=\ve(1-\ve/2+\ve^2/3)+\O(\ve^4)
\]
(\ref{eq:C1-first-asymp}) and (\ref{eq:eta-defined}):
\begin{equation}\label{eq:eta-asymp}
\eta=\frac{\gamma}{2} +\frac{\ve}{2\pi}\ln\frac{1}{\ve}+\O(\ve).
\end{equation}
For theta functions we get from
(\ref{eq:theta1-def})-(\ref{eq:theta3-def}) for $q_1\to{0}$:
\begin{equation}\label{eq:theta1-asymp}
\vt_1(\eta,q_1)=2q_1^{1/4}(\sin(\eta)+\O(q_1^{2})),
\end{equation}
\begin{equation}\label{eq:theta2-asymp}
\vt_2(\eta,q_1)=2q_1^{1/4}(\cos(\eta)+\O(q_1^{2}))
\end{equation}
\begin{equation}\label{eq:theta3-asymp}
\vt_3(\eta,q_1)=1+2q_1\cos(2\eta)+\O(q_1^4),
\end{equation}
\begin{equation}\label{eq:theta4-asymp}
\vt_4(\eta,q_1)=1-2q_1\cos(2\eta)+\O(q_1^4),
\end{equation}
\begin{equation}\label{eq:theta2prime-asymp}
\vt_2'(\eta,q_1)=2q_1^{1/4}(-\sin(\eta)+\O(q_1^{2})).
\end{equation}
For the complete elliptic integral $K$ we can write according to
(\ref{eq:omega1-fixed}) and (\ref{eq:theta3-newnome}) as
$q_1\to{0}$:
\begin{equation}\label{eq:K-asymp}
K(q_1)=\frac{\pi}{2}\vt_3^2(0,q)=\frac{1}{2}\ln(1/q_1)\vt_3^2(0,q_1)=\frac{1}{2}\ln\frac{1}{q_1}\left(1+2\sum\limits_{n=1}^{\infty}q_1^{n^2}\right)^2=
\frac{1}{2}\ln(1/q_1)(1+4q_1+\O(q_1^2)).
\end{equation}

For moduli $k$ and $k'=\sqrt{1-k^2}$ we get according to
(\ref{eq:k-defined}), (\ref{eq:k'-defined}),
(\ref{eq:theta2-newnome}), (\ref{eq:theta3-newnome}) and
(\ref{eq:theta2-asymp})-(\ref{eq:theta4-asymp}):
\[
k=\frac{\vt_4^2(0,q_1)}{\vt_3^2(0,q_1)}=\frac{(1-2q_1+\O(q_1^4))^2}{(1+2q_1+\O(q_1^4))^2}=
\frac{1-4q_1+4q_1^2+\O(q_1^4)}{1+4q_1+4q_1^2+\O(q_1^4)}=
\]\[
=(1-4q_1+4q_1^2+\O(q_1^4))(1-4q_1-4q_1^2+16(q_1+q_1^2+\O(q_1^4))^2-64(q_1+\O(q_1^2))^3+\O(q_1^4))=
\]
\begin{equation}\label{eq:k-asymp}
=(1-4q_1+4q_1^2+\O(q_1^4))(1-4q_1+12q_1^2-32q_1^3+\O(q_1^4))=1-8q_1+32q_1^2-96q_1^3+\O(q_1^4).
\end{equation}
Hence for $q_1\to{0}$:
\begin{equation}\label{eq:k2-asymp}
k^2=1-16q_1+128q_1^2-704q_1^3+\O(q_1^4),
\end{equation}
\begin{equation}\label{eq:sqrtk-asymp}
\sqrt{k}=1-4q_1+8q_1^2-16q_1^3+\O(q_1^4),
\end{equation}
and
\begin{equation}\label{eq:recsqrtk-asymp}
1/\sqrt{k}=1+4q_1+8q_1^2+16q_1^3+\O(q_1^4),
\end{equation}
where the expansions (valid for $z\to{0}$):
\begin{equation}\label{eq:expansions}
(1+z)^r=1+rz+\binom{r}{2}z^2+\binom{r}{3}z^3+\binom{r}{4}z^4+\O(z^5),~~~\frac{1}{1-z}=1+z+z^2+z^3+z^4+\O(z^5).
\end{equation}
have been used.

In a similar fashion:
\begin{equation}\label{eq:k'-asymp}
k'=\frac{\vt_2^2(0,q_1)}{\vt_3^2(0,q_1)}=4q_1^{1/2}(1-4q_1+14q_1^2-40q_1^3+\O(q_1^4)),
\end{equation}
\begin{equation}\label{eq:sqrtk'-asymp}
\sqrt{k'}=2q_1^{1/4}(1-2q_1+5q_1^2-10q_1^3+\O(q_1^4))
\end{equation}
and
\begin{equation}\label{eq:sqrtk'k-asymp}
\sqrt{k'/k}=2q_1^{1/4}(1+2q_1+5q_1^2+10q_1^3+\O(q_1^4)).
\end{equation}
Then from (\ref{eq:sn-q1})-(\ref{eq:Z-q1}) and
(\ref{eq:theta1-asymp})-(\ref{eq:K-asymp}):
\begin{equation}\label{eq:sn-asymp}
\sn(i\alpha,k)=i\tan(\eta)+4iq_1\tan(\eta)+\O(q_1^{2}),
\end{equation}
\begin{equation}\label{eq:ksn-asymp}
-k^2\sn^2(i\alpha,k)=\tan^2(\eta)-8q_1\tan^2(\eta)+\O(q_1^2),
\end{equation}
\begin{equation}\label{eq:cn-asymp}
\cn(i\alpha,k)=\frac{1}{\cos(\eta)}+4q_1\frac{\sin^2(\eta)}{\cos(\eta)}+\O(q_1^2),
\end{equation}
\begin{equation}\label{eq:dn-asymp}
\dn(i\alpha,k)=\frac{1}{\cos(\eta)}-4q_1\frac{\sin^2(\eta)}{\cos(\eta)}+\O(q_1^2),
\end{equation}
\begin{equation}\label{eq:Z-asymp}
\Z(i\alpha,k)=i\tan(\eta)-\frac{2i\eta}{\ln(1/q_1)}-4iq_1\tan(\eta)+\frac{8iq_1\eta}{\ln(1/q_1)}+\O(q_1^2),
\end{equation}
\begin{equation}\label{eq:cndn-sn-asymp}
2i\frac{\cn(i\alpha)\dn(i\alpha)}{\sn(i\alpha)}=\frac{2}{\sin(\eta)\cos(\eta)}-\frac{8q_1}{\sin(\eta)\cos(\eta)}+\O(q_1^2).
\end{equation}
Introduce the notation
\begin{equation}\label{eq:x-defined}
x=[\ln(1/q_1)]^{-1}=\frac{1}{\pi\cp(C_1')}\to{0}~~\text{as}~q_1\to{0}.
\end{equation}
Using this notation and (\ref{eq:lambda-defined}):
\[
\lambda^2(R,q_1)=\frac{(i\tan(\eta)-2i\eta{x}+\O(q_1))/[(1+\O(q_1))(i\tan(\eta)+\O(q_1))]}
{\left(\frac{1}{\cos(\eta)}+\O(q_1)\right)\left(\frac{1}{\cos(\eta)}+\O(q_1)\right)+(i\tan(\eta)+\O(q_1))(i\tan(\eta)-2i\eta{x}+\O(q_1))}
\]\[
=\frac{1-2\eta{x}/\tan\eta+\O(q_1)}
{1/\cos^2\eta-\tan^2\eta+2\eta{x}\tan\eta+\O(q_1)}
\]
\begin{equation}\label{eq:lambda2-asymp}
=\frac{1-2\eta{x}/\tan\eta}{1+2\eta{x}\tan\eta}+\O(q_1)=1-\frac{2\eta{x}}{\sin\eta(\cos\eta+2\eta{x}\sin\eta)}+\O(q_1).
\end{equation}
By (\ref{eq:expansions}) this leads to:
\[
\lambda(R,q_1)=1+\sum\limits_{k=1}^{\infty}(-1)^k\binom{1/2}{k}
\frac{(2\eta{x})^k}{\sin^k\eta\cos^k\eta(1+2\eta{x}\tan\eta)^k}+\O(q_1)
\]\[
=\sum\limits_{k=0}^{\infty}\binom{1/2}{k}
\frac{(-2\eta{x}\tan\eta)^k}{\sin^{2k}\eta}\sum\limits_{m=0}^{\infty}\binom{k+m-1}{m}(-2\eta{x}\tan\eta)^m+\O(q_1)
\]\[
=\sum\limits_{n=0}^{\infty}(-2\eta{x}\tan\eta)^n\sum\limits_{k=0}^{n}\binom{1/2}{k}\binom{n-1}{n-k}\sin^{-2k}\eta
+\O(q_1),
\]
where $\binom{n-1}{n}=0$ when $n>0$ and $\binom{n-1}{n}=1$, when
$n=0$ (which is standard and very convenient convention). From
this:
\[
\sigma\equiv1-\lambda(R,q_1)=-\sum\limits_{n=1}^{\infty}(-2\eta{x}\tan\eta)^n\sum\limits_{k=1}^{n}\binom{1/2}{k}\binom{n-1}{n-k}\sin^{-2k}\eta
+\O(q_1)
\]
\begin{equation}\label{eq:delta-defined}
=\frac{\eta{x}\tan\eta}{\sin^2\eta}(1-\beta(x,\eta))+\O(q_1),
\end{equation}
where
\[
\beta(x,\eta)=4\eta{x}\tan\eta\sum\limits_{n=0}^{\infty}(-1)^n(2\eta{x}\tan\eta)^n
\sum\limits_{k=0}^{n+1}\binom{1/2}{k+1}\binom{n+1}{n-k+1}\sin^{-2k}\eta=2\eta{x}\tan\eta\left(1-\frac{1}{4\sin^{2}\eta}\right)
\]\[
-(2\eta{x}\tan\eta)^2\left(1-\frac{1}{2\sin^{2}\eta}+\frac{1}{8\sin^{4}\eta}\right)+
(2\eta{x}\tan\eta)^3\left(1-\frac{3}{4\sin^{2}\eta}+\frac{3}{8\sin^{4}\eta}-\frac{5}{64\sin^{6}\eta}\right)
\]
\begin{equation}\label{eq:beta-defined}
-(2\eta{x}\tan\eta)^4\left(1-\frac{1}{\sin^{2}\eta}+\frac{3}{4\sin^{4}\eta}-\frac{5}{16\sin^{6}\eta}+
\frac{7}{128\sin^{8}\eta}\right)+ \O(x^5).
\end{equation}
 Now denote:
\begin{equation}\label{eq:nu-defined}
\nu\equiv-k^2\sn^2(i\alpha,k)=\tan^2(\eta)(1-8q_1+\O(q_1^2))>0~~\text{when}~q_1~\text{is~small},
\end{equation}
\begin{equation}\label{eq:sqrtnu-asymp}
\sqrt{\nu}=\tan(\eta)(1-4q_1+\O(q_1^2)).
\end{equation}
We have by (\ref{eq:k-asymp}) and (\ref{eq:delta-defined}):
\[
(1-k)/\sigma=\O(q_1\ln(1/q_1))\to{0}~\text{as}~q_1\to{0}.
\]
It follows from the results of \cite{KSS} that for
\begin{equation}\label{eq:k-faster}
k\to{1},~~~\sigma\to{0}~~\text{and}~~(1-k)/\sigma\to{0}
\end{equation}
$\Pi-F$ has the asymptotic approximation
\begin{equation}\label{eq:Pi-Fasymp}
\Pi(1-\sigma,\nu,k)-F(1-\sigma,k)=\frac{-\nu}{2(1+\nu)}\ln\frac{2-\sigma}{\sigma}+
\frac{\sqrt{\nu}\arctan\left((1-\sigma)\sqrt{\nu}\right)}{1+\nu}+\O\left((1-k)/\sigma\right).
\end{equation}
The Taylor expansion for $\arctan(\tan\eta(1-\epsilon))$ for small
$\epsilon$ is given by
\[
\arctan(\tan\eta(1-\epsilon))=\eta+\sum\limits_{n=1}^{\infty}\frac{1}{n!}\frac{d^n}{dx^n}\arctan(x)_{|x=\tan\eta}(-\epsilon\tan\eta)^n.
\]
Simple manipulations reveal
\[
\frac{d^n}{dx^n}\arctan(x)_{|x=\tan\eta}=-(n-1)!\cos^{2n}\eta\sum\limits_{m=0}^{[(n-1)/2]}(-1)^{n-m}\binom{n}{2m+1}\tan^{n-2m-1}\eta.
\]
Hence
\[
\arctan(\tan\eta(1-\epsilon))=\eta-\sum\limits_{n=1}^{\infty}\frac{\epsilon^n}{n}\sin^n\eta\sum\limits_{m=0}^{[(n-1)/2]}(-1)^{m}\binom{n}{2m+1}\sin^{n-2m-1}\eta\cos^{2m+1}\eta
\]\[
=\eta-\epsilon\cos\eta\sin\eta-
\epsilon^2\cos\eta\sin^3\eta-\frac{\epsilon^3}{3}\cos\eta\sin^3\eta(4\sin^2\eta-1)
-\epsilon^4\cos\eta\sin^5\eta(2\sin^2\eta-1)+ \O(\epsilon^5).
\]
The Taylor expansion for $\ln(2-\sigma)$ is:
\[
\ln(2-\sigma)=\ln(2)-\sum\limits_{m=1}^{\infty}\frac{\sigma^m}{m2^m}.
\]
Using these Taylor expansions, (\ref{eq:nu-defined})
 and (\ref{eq:sqrtnu-asymp}) expansion (\ref{eq:Pi-Fasymp})
is transformed into
\begin{multline*}
\Pi(\lambda,\nu,k)-F(\lambda,k)=
\frac{-\tan^2\eta(\ln(1/\sigma)+\ln(2))}{2(1+\tan^2\eta+\O(q_1))}
+\frac{(\tan\eta+\O(q_1))\arctan\left((1-\sigma)(\tan\eta+\O(q_1))\right)}{1+\tan^2\eta+\O(q_1)}\\[0pt]
+\frac{\tan^2\eta+\O(q_1)}{1+\tan^2\eta+\O(q_1)}\left(\sum\limits_{m=1}^{\infty}\frac{\sigma^m}{m2^{m+1}}\right)+\O\left((1-k)/\sigma\right)
=-\frac{1}{2}\sin^2\eta\ln(1/\sigma)-\frac{1}{2}\sin^2\eta\ln2\\[0pt]
+\frac{\tan\eta\arctan\left(\tan\eta[1-\sigma+\O(q_1)]\right)}{1+\tan^2\eta}+
\sin^2\eta\left(\sum\limits_{m=1}^{\infty}\frac{\sigma^m}{m2^{m+1}}\right)+\O(q_1\ln(1/q_1))=
-\frac{1}{2}\sin^2\eta\ln(1/\sigma)\\[0pt]
-\frac{\sin^2\eta\ln2}{2}+\frac{\eta\sin(2\eta)}{2}
+\sum\limits_{n=1}^{\infty}\frac{\sigma^n}{n}\!
\left[\frac{\sin^2\eta}{2^{n+1}}-\!\!\!\sum\limits_{m=0}^{[(n-1)/2]}\!\!(-1)^{m}\binom{n}{2m+1}\cos^{2m+2}\eta\sin^{2n-2m}\eta\right]
+\O\!\left(q_1\ln\frac{1}{q_1}\!\right)\\[0pt]
=-\frac{1}{2}\sin^2\eta\ln\frac{1}{\sigma}-\frac{1}{2}\sin^2\eta\ln(2)+\frac{1}{2}\eta\sin(2\eta)+
\sigma\sin^2\eta(1/4-\cos^2\eta)+\sigma^2\sin^2\eta(1/16-\cos^2\eta\sin^2\eta)\\[5pt]
+\sigma^3\sin^2\eta(1/48-\cos^2\eta\sin^4\eta+\cos^4\eta\sin^2\eta/3)
+\sigma^4\sin^2\eta(1/128-\cos^2\eta\sin^6\eta+\cos^4\eta\sin^4\eta)+\O\left(\sigma^5\right).
\end{multline*}
This and (\ref{eq:cndn-sn-asymp}) lead to
\begin{multline}\label{eq:cndnsnPiF}
2i\frac{\cn(i\alpha)\dn(i\alpha)}{\sn(i\alpha)}(\Pi(1-\sigma,\nu,k)-F(1-\sigma,k))
=-\tan\eta\ln\frac{1}{\sigma}-\tan\eta\ln(2)+2\eta+2\sigma\tan\eta(1/4-\cos^2\eta)\\
+2\sigma^2\tan\eta(1/16-\cos^2\eta\sin^2\eta)+2\sigma^3\tan\eta(1/48-\cos^2\eta\sin^4\eta+\cos^4\eta\sin^2\eta/3)\\
+2\sigma^4\tan\eta(1/128-\cos^2\eta\sin^6\eta+\cos^4\eta\sin^4\eta)+\O\left(\sigma^5\right).
\end{multline}
Substituting (\ref{eq:delta-defined}) for $\sigma$ and denoting
for brevity
\[
z=2\eta{x}\tan\eta,~~~~s=\frac{1}{\sin^2\eta},
\]
 we obtain:
\[
\ln\frac{1}{\sigma}=\ln\frac{1}{x}+\ln\frac{\sin2\eta}{2\eta}+\sum\limits_{k=1}^{\infty}\frac{\beta^k(x,\eta)}{k}=
\]\[
=\ln\frac{1}{x}+\ln\frac{\sin2\eta}{2\eta}+z(1-s/4)+z^2(-1/2+s/4-3s^2/32)+z^3(1/3-s/4+3s^2/16-5s^3/96)
\]\[
+z^4(-1/4+s/4-9s^2/32+5s^3/32-35s^4/1024)+\O(x^5).
\]
For the powers of $\sigma$ we compute  up to $\O(x^5)$:
\[
\sigma=zs/2-z^2(s/2-s^2/8)+z^3(s/2-s^2/4+s^3/16)-z^4(s/2-3s^2/8+3s^3/16-5s^4/128)+\O(x^5),
\]
\[
\sigma^2=z^2s^2/4-z^3(s^2/2-s^3/8)+z^4(3s^2/4-3s^3/8+5s^3/64)+\O(x^5),
\]
\[
\sigma^3=z^3s^3/8-z^4(3s^3/8-3s^4/32)+\O(x^5),
\]
\[
\sigma^4=z^4s^4/16+\O(x^5).
\]
Hence (\ref{eq:cndnsnPiF})  transforms into
\begin{multline*}
2i\frac{\cn(i\alpha)\dn(i\alpha)}{\sn(i\alpha)}(\Pi(1-\sigma,\nu,k)-F(1-\sigma,k))=
-\tan\eta\ln\frac{1}{x}+2\eta-\tan\eta\ln\frac{\sin2\eta}{\eta}
\\
-\frac{\eta{x}}{\cos^2{\eta}}-\frac{\eta^2x^2(1-4\sin^2\eta)}{4\cos^3\eta\sin\eta}
-\frac{\eta^3x^3(1-4\sin^2\eta+8\sin^4\eta)}{6\cos^4\eta\sin^2\eta}\\
-\frac{\eta^4x^4(10-109\sin^2\eta+208\sin^4\eta-144\sin^6\eta)}{32\sin^3\eta\cos^5\eta}+\O(x^5).
\end{multline*}
For $2iF(\lambda,k)Z(i\alpha,k)$ from (\ref{eq:Z-asymp}) and
\[
F(1-\sigma,k)=\frac{1}{2}\ln(1/\sigma)+\frac{1}{2}\ln(2)-\sum\limits_{m=1}^{\infty}\frac{\sigma^m}{m2^{m+1}}+\O(q_1\ln(q_1))
\]
we get representation
\begin{multline*}
2iF(\lambda,k)Z(i\alpha,k)=-\tan\eta\ln(1/x)-\tan\eta\ln\frac{\sin(2\eta)}{\eta}+2\eta{x}\ln(1/x)
+\eta{x}\left(\frac{\cos(2\eta)}{\cos^2\eta}+2\ln\frac{\sin(2\eta)}{\eta}\right)
\\
-\frac{\eta^2x^2(5-16\sin^2\eta+8\sin^4\eta)}{4\cos^3\eta\sin\eta}-\frac{\eta^3x^3(4-15\sin^2\eta+24\sin^4\eta-8\sin^6\eta)}{6\cos^4\eta\sin^2\eta}
\\
-\frac{\eta^4x^4(55-256\sin^2\eta+480\sin^4\eta-512\sin^6\eta+128\sin^8\eta)}{96\sin^3\eta\cos^5\eta}+\O(x^5).
\end{multline*}
Substituting last two expansions into (\ref{eq:gamma*-found}) we
finally arrive at the following formula:
\begin{multline}\label{eq:gamma-asymp}
\gamma=2\eta-2\eta{x}\ln(1/x)-2\eta{x}\left(1+\ln\frac{\sin(2\eta)}{\eta}\right)+2\eta^2x^2\cot(2\eta)
\\
+2\eta^3x^3\frac{3-2\sin^2(2\eta)}{3\sin^2(2\eta)}+\eta^4x^4\frac{(25+96\sin^2\eta-48\sin^4\eta-128\sin^6\eta)}{12\sin^3(2\eta)}+\O(x^5).
\end{multline}

Now denote $y=\cp(C_1')$ and recall the rough approximations
(\ref{eq:C1-first-asymp}) and (\ref{eq:eta-asymp}). These
asymptotic approximations and equation (\ref{eq:gamma-asymp}) are
sufficient ingredients to start using the bootstrapping technique
for identifying further asymptotic terms in
(\ref{eq:C1-first-asymp}). Substituting (\ref{eq:C1-first-asymp})
and (\ref{eq:eta-asymp}) into (\ref{eq:gamma-asymp}) we can
rewrite the latter  as
\begin{equation}\label{eq:gamma-asymp-y}
\gamma=\ln(1+\ve)\left[y-\frac{1}{\pi}\ln(\pi)-\frac{1}{\pi}\ln{y}-\frac{1}{\pi}
\left(1+\ln\frac{2\sin(\gamma+\frac{\ve}{\pi}\ln\frac{1}{\ve}+\O(\ve))}{\gamma+\frac{\ve}{\pi}\ln\frac{1}{\ve}+\O(\ve)}\right)\right]
+\O(\ve^2).
\end{equation}
The following Taylor expansions hold true:
\[
\ln{y}=\ln\left[\frac{\gamma}{\ve}\left(1+\frac{\ve}{\gamma\pi}\ln\frac{1}{\ve}+\O(\ve)\right)\right]=
\ln\frac{\gamma}{\ve}+\frac{\ve}{\gamma\pi}\ln\frac{1}{\ve}+\O(\ve),
\]
\[
\frac{\sin(\gamma+\frac{\ve}{\pi}\ln\frac{1}{\ve}+\O(\ve))}{\gamma+\frac{\ve}{\pi}\ln\frac{1}{\ve}+\O(\ve)}=
\frac{\sin\gamma}{\gamma}+\frac{\ve}{\pi}\left(\frac{\cos\gamma}{\gamma}-\frac{\sin\gamma}{{\gamma}^2}\right)\ln\frac{1}{\ve}+\O(\ve),
\]
\[
\ln\frac{2\sin(\gamma+\frac{\ve}{\pi}\ln\frac{1}{\ve}+\O(\ve))}{\gamma+\frac{\ve}{\pi}\ln\frac{1}{\ve}+\O(\ve)}=
\ln\left[\frac{2\sin\gamma}{\gamma}\left(1+\frac{\ve\gamma}{\pi\sin\gamma}
\left(\frac{\cos\gamma}{\gamma}-\frac{\sin\gamma}{{\gamma}^2}\right)\ln\frac{1}{\ve}+\O(\ve)\right)\right]=
\]\[
=\ln\frac{2\sin\gamma}{\gamma}
+\frac{\ve}{\pi}\left(\cot\gamma-\frac{1}{\gamma}\right)\ln\frac{1}{\ve}+\O(\ve).
\]
Substituting these expansions into (\ref{eq:gamma-asymp-y})
yields:
\[
\frac{\gamma}{\ve}=\left(1-\frac{\ve}{2}\right)
\left(y-\frac{1}{\pi}\ln\frac{1}{\ve}
-\frac{1}{\pi}(1+\ln(2\pi\sin\gamma))
-\frac{\cot\gamma}{\pi^2}\ve\ln\frac{1}{\ve}\right) +\O(\ve).
\]
Substituting once more $\gamma+\ve\ln(1/\ve)/\pi+\O(\ve)$ for
$y\ve$ results in identification of the next asymptotic terms:
\begin{equation}\label{y-second}
y=\frac{\gamma}{\ve}+\frac{1}{\pi}\ln\frac{1}{\ve}+\frac{\gamma}{2}+\frac{1}{\pi}
+\frac{1}{\pi}\ln(2\pi\sin\gamma)+\frac{\cot\gamma}{\pi^2}\ve\ln\frac{1}{\ve}+\O(\ve)
\end{equation}
and
\begin{equation}\label{eta-second}
\eta=\frac{\gamma}{2}+\frac{\ve}{2\pi}\ln\frac{1}{\ve}+\frac{\ve}{2\pi}(1+\ln(2\pi\sin\gamma))
+\frac{\ve^2}{2\pi}\left(\frac{1}{\pi}\cot\gamma-\frac{1}{2}\right)\ln\frac{1}{\ve}+\O(\ve^2).
\end{equation}
The next step of bootstrapping is to substitute these
approximations into (\ref{eq:gamma-asymp}).  Using the expansions
(as $x\to{0}$):
\[
\ln(1+x)=x(1-x/2+x^2/3)+\O(x^4),
\]\[
\frac{\sin(\gamma+x)}{\gamma+x}=\frac{\sin\gamma}{\gamma}
+\frac{x}{\gamma}\left(\cos\gamma-\frac{\sin\gamma}{\gamma}\right)
-\frac{x^2}{{\gamma}^2}\left(\frac{\gamma}{2}\sin\gamma+\cos\gamma-\frac{\sin\gamma}{\gamma}\right)+\O(x^3),
\]
\[
\cot(\gamma+x)=\cot\gamma-\frac{x}{\sin^2\gamma}+\O(x^2)
\]
and omitting some tedious calculations we arrive at:
\[
y=\frac{\gamma}{\ve}+\frac{1}{\pi}\ln\frac{1}{\ve}+\frac{\gamma}{2}+\frac{1}{\pi}(1
+\ln(2\pi\sin\gamma))+\frac{\cot\gamma}{\pi^2}\ve\ln\frac{1}{\ve}
+\ve\left[\frac{\gamma}{4}+\frac{1}{2\pi}+\frac{\cot\gamma}{\pi^2}\left(\frac{1}{2}+\ln(2\pi\sin\gamma)\right)\right]
\]
\begin{equation}\label{eq:y-third}
-\frac{\ve^2\ln^2(1/\ve)}{2\pi^3\sin^2\gamma}
-\frac{\ve^2\ln(1/\ve)}{2\pi^3}\left(1+\pi\cot\gamma-\cot^2\gamma+\frac{2\ln(2\pi\sin\gamma)}{\sin^2\gamma}\right)
+\O(\ve^2).
\end{equation}
To derive an asymptotic expansion for the capacity of the
condenser $C_1$ (see section~2) from this formula recall that
$\cp(C_1)=\cp(C_1')$ and substitute (\ref{eq:eps-delta}) for
$\ve$.  Simple calculations lead to the expression
\[
\cp(C_1)=\frac{\gamma}{\delta}+\frac{1}{\pi}\ln\frac{1}{\delta}+\frac{1}{\pi}(1
+\ln(2\pi\sin\gamma))+\frac{\cot\gamma}{\pi^2}\delta\ln\frac{1}{\delta}
+\frac{\delta\cot\gamma}{\pi^2}\left(\frac{1}{2}+\ln(2\pi\sin\gamma)\right)
\]
\begin{equation}\label{eq:C-1-fullasymp}
-\frac{\delta^2\ln^2(1/\delta)}{2\pi^3\sin^2\gamma}
-\frac{\delta^2\ln(1/\delta)}{2\pi^3\sin^2\gamma}\left(2\ln(2\pi\sin\gamma)-\cos(2\gamma)\right)
+\O(\delta^2).
\end{equation}
The final step that gives formula (\ref{eq:C-rho-asymp}) from
Theorem~\ref{th:main} is to carry out the conformal mapping
$z\to\rho{z}$. Under this mapping the condenser $C_1$ transforms
into the condenser $C_\rho$ with $L=2\rho\gamma$ and
$h=2\delta\rho$. Thus substituting $\gamma=L/(2\rho)$ and
$\delta=h/(2\rho)$ into (\ref{eq:C-1-fullasymp}) brings the
desired result.

\paragraph{5. Acknowledgements.}  The author thanks Professor
V.N.~Dubinin for introducing him into the subject of capacities
and for general guidance.  This research has been supported by the
Russian Science Support Foundation, the Far Eastern Branch of the
Russian Academy of Sciences (grant 04-3-$\Gamma$-01-020) and
Russian Basic Research Fund (grant 02-01-000-28).

\end{document}